# CLOSURE OF RIGID SEMIANALYTIC SETS

Hans Schoutens

ABSTRACT. Let $K$ be an algebraically closed field of characteristic zero, endowed with a complete non-archimedean norm. Let $X$ be a $K$-rigid analytic variety and $\Sigma$ a semianalytic subset of $X$. Then the closure of $\Sigma$ in $X$ with respect to the canonical topology is again semianalytic. The proof uses Embedded Resolution of Singularities.

## 0. Introduction

Fix an algebraically closed field $K$ of characteristic zero which is endowed with a complete non-archimedean norm. Over $K$ are defined the so called *rigid analytic varieties*; see [**BGR**] for their definition and various properties. They locally look like an *affinoid variety* $\operatorname{Sp} A$, where $A$ is an *affinoid algebra*, that is, a quotient of a free Tate ring $K\langle S\rangle$. The latter algebra is defined as the collection of all power series over $K$ in the variables $S = (S_1, \dots, S_n)$ for which the general coefficient tends to zero. If $R$ denotes the valuation ring of $K$, that is, all elements of norm at most 1, then $K\langle S\rangle$ consists exactly of those power series converging on $R^n$.

There is a natural topology on the spaces $R^n$ induced by the norm, called the *canonical* topology. To be more explicit, for an $n$-tuple $x = (x_1, \dots, x_n) \in R^n$, we define

$$|x| = \max_{1 \leq i \leq n} \{|x_i|\}.$$

Any affinoid space $X = \operatorname{Sp} A$ can be considered as a closed subvariety of some $R^n$ and hence is induced with a canonical topology (which is independent from the particular choice of the closed immersion, see [**BGR**, 7.2.1] ). Nevertheless, in order to perform a good analytic theory on these spaces, one needs to work with a coarser topology, since affinoid varieties are totally disconnected in the canonical topology. This topology will no longer be a classical one, but only a Grothendieck topology. We refer to [**BGR**] for more details. Since any admissible open in this Grothendieck topology is also open in the canonical topology, we can also define the canonical topology on an arbitrary rigid analytic variety $X$, by saying that a subset $U \subset X$ is open in the canonical topology, if and only if, there exists an admissible affinoid covering $\{X_i\}_i$ of $X$, such that each $U \cap X_i$ is open in the canonical topology on $X_i$. For an arbitrary subset $\Sigma$ of a rigid analytic variety $X$, we denote by $\operatorname{cl}(\Sigma)$ its closure with respect to the canonical topology.

In this paper, we study the closure in the canonical topology of a semianalytic subset of a rigid analytic variety. Let us briefly recall the definition of a semianalytic set and introduce some notation. First, let $X = \operatorname{Sp} A$ be an affinoid variety. We define the class $\mathcal{S}_X$ of all subsets of the form

$$(1) \qquad \Sigma = \{\, x \in X \mid |f(x)| \diamond |g(x)| \,\},$$

where $f, g \in A$ and $\diamond \in \{\leq, <, \geq, >\}$. By restricting to those sets for which either $f$ or $g$ is a unit in the description (1) of $\Sigma$, we obtain a subclass, denoted $\mathcal{A}_X$. Note that without loss of generality we then always may put this unit equal to 1. We call $\Sigma$ a *basic set* of $X$, if it is a finite intersection of sets in $\mathcal{S}_X$. We call $\Sigma$ *globally semianalytic*, if it is a finite union of basic subsets. The class of globally semianalytic subsets of $X$ forms a Boolean algebra, generated by the subclass $\mathcal{S}_X$.

There exist also some more local versions. First, let $\mathcal{S}_X^{\operatorname{loc}}$ (respectively, $\mathcal{A}_X^{\operatorname{loc}}$) denote the collection of all subsets $\Sigma$ of $X$ for which there exists a finite admissible affinoid covering $\{X_i\}_i$ of $X$, such that $\Sigma \cap X_i$ lies in $\mathcal{S}_{X_i}$ (respectively, in $\mathcal{A}_{X_i}$), for every $i$. The set $\Sigma$ is called *semianalytic*, if it is a finite union of finite intersections of sets out of $\mathcal{S}_X^{\operatorname{loc}}$.

---

For an arbitrary rigid analytic variety $X$, we define similarly $\mathcal{S}_X^{\text{loc}}$ (respectively, $\mathcal{A}_X^{\text{loc}}$) as the collection of all subsets $\Sigma$ of $X$, for which there exists an admissible affinoid covering $\{X_i\}_i$ of $X$, such that, for each $i$, we have that $\Sigma \cap X_i$ belongs to $\mathcal{S}_{X_i}$ (respectively, belongs to $\mathcal{A}_{X_i}$). Semianalytic sets are likewise defined.

One defines a *subanalytic* set as a projection of a semianalytic set. A special subclass, termed *strongly* subanalytic sets are defined and studied in [**Sch 1**] – [**Sch 3**], whereas the general case is treated in [**GS**]. The following now holds.

**Theorem.** *Let $X$ be a rigid analytic variety and $\Sigma$ a (strongly) subanalytic subset of $X$. Then its closure $\text{cl}(\Sigma)$ in the canonical topology is also (strongly) subanalytic.*

*Proof.* For the strong version see [**Sch 1**, Corollary 5.4] and for the general case see [**LR 2**, Corollary 1.2]. ∎

The question remained, whether the same holds true, when we replace *(strongly) subanalytic* by *semianalytic*. This is the contents of our main theorem **(3.3)**.

**Main Theorem.** *Let $X$ be a rigid analytic variety and $\Sigma$ a semianalytic subset of $X$. Then its closure $\text{cl}(\Sigma)$ in the canonical topology is also semianalytic.*

An analogous theorem for real semianalytic sets had been proven in [**Łoj 2**]. A $p$-adic version, has recently been shown by LIU in his preprint [**Liu**], using cell decomposition. After submission, the author learned that HUBER has proven the same theorem (without restriction on the characteristic) in [**Hub**]. Our proof uses entirely different methods and it could be easily adapted to the $p$-adic situation to obtain LIU's result. We do make use, however, of a corollary to the highly non-trivial theorem of Embedded Resolution of Singularities [**Sch 5**, Theorem 3.2.3], which is proved using HIRONAKA's general Embedded Resolution. The drawback of using this method, is that we have to confine ourselves to the zero characteristic case, since this theorem is not (yet?) known in positive characteristic.

In the first section, we study more closely blowing up maps, which appear in the resolution of singularities. For their definition and elementary properties in the rigid analytic setting, we refer to our article [**Sch 4**]. Here we prove that blowing up maps are algebraic and they send subanalytic sets which are closed for the canonical topology to (subanalytic) sets which are closed for the canonical topology. The second section is concerned with the closure of a *constructible* set (i.e., a finite union of differences of two closed analytic subsets). Recall that a subset $V$ of a rigid analytic variety $X$ is called a *closed analytic subset*, if their exists a coherent $\mathcal{O}_X$-ideal $\mathcal{I}$, such that $V$ consists exactly of those points for which $\mathcal{I}_x \neq \mathcal{O}_{X,x}$. It is shown that the closure of a constructible set in the canonical topology is closed analytic. In the third section, we then provide a proof for our Main Theorem, using these previous results together with the already mentioned Embedded Resolution of Singularities.

In the last section, we investigate a little bit closer the boundary of some semianalytic sets. With the *boundary* $\partial \Sigma$ (with respect to the canonical topology) of a set $\Sigma$, we mean the difference of the closure $\text{cl}(\Sigma)$ (in the canonical topology) with its interior (in the canonical topology). We obtain the following remarkable fact.

**Boundary Theorem.** *Let $X$ be a rigid analytic variety. Let*

$$\Sigma = \{\, x \in R^n \mid |f(x)| < |g(x)| \,\},$$

*where $f, g \in K\langle S \rangle$. If $f$ and $g$ are relatively prime, then the boundary of $\Sigma$ is the analytic subset $V(f, g)$.*

*Acknowledgment.* We would like to thank JAN DENEF, who suggested the use of Embedded Resolution of Singularities as a tool to prove our Main Theorem and the referee for pointing out a mistake in Theorem **(4.4)**.

**0.2. Preliminaries.** Recall that a subset $A$ in a topological space $X$ is called *clopen*, if it is both open and closed. It is easy to see that the closure of a finite union of sets is the union of the closures of these sets. The same statement for intersection instead of union is no longer true. However, if $A$ is clopen and $B$ arbitrary, then it easy to verify that

$$(0.1) \qquad \text{cl}(A \cap B) = A \cap \text{cl}(B).$$

Let $h\colon X \to Y$ be a continuous map of topological spaces and let $\Sigma \subset Y$. Then a straightforward argument shows that

$$h(\operatorname{cl}(h^{-1}(\Sigma))) \subset \operatorname{cl}(\Sigma). \tag{0.2}$$

We will also be making use of the following theorem due to LIPSHITZ.

**0.3. Proposition.** *Let $\Gamma$ and $\Sigma$ be subanalytic subsets of $R^n$ which are closed for the canonical topology. If $\Gamma \cap \Sigma = \emptyset$, then there exists a $\delta > 0$, such that for all $q \in \Gamma$ and $s \in \Sigma$, we have that $|q - s| \geq \delta$. In other words, the distance between $\Gamma$ and $\Sigma$ is strictly positive.*

*Proof.* This is merely a restatement of [**Lip**, Corollary 5.9]. Consider the following function $f\colon \Gamma \to \mathbb{R}^+$ given by $q \mapsto d(q, \Sigma)$, where $d(q, \Sigma)$ stands for the distance from $q$ to $\Sigma$, given as the infimum of all $|q - s|$, when $s$ runs over $\Sigma$. Since $\Sigma$ is closed in the canonical topology, one easily checks that $f$ nowhere vanishes on $\Gamma$. It is also an easy exercise to show that $f$ is a subanalytic and continuous function. We therefore are entitled to apply [**Lip**, Corollary 5.9] to conclude that

$$d(\Gamma, \Sigma) = \inf_{q \in \Gamma} f(q) \neq 0.$$

(Note that LIPSHITZ's notion of subanalyticity is broader than ours). ∎

## 1. BLOWING UP MAPS

**1.1. Lemma.** *A rigid analytic blowing up map is algebraic.*

*Proof.* Saying that a map $h : X \to Y$ is algebraic amounts to saying that its graph $\Gamma(h)$ is an algebraic subset of $X \times Y$. So let $h : \tilde{X} \to X$ be the blowing up of $X$ with centre $Z$. Without loss of generality, we can assume that $X$ is affinoid. By [**Sch 4**, 2.2.2], we can find a closed immersion $X \hookrightarrow R^n \times R^k$, such that $Z = X \cap (0 \times R^k)$. If $\pi : V \to R^n \times R^k$ is the blowing up with centre $0 \times R^k$, then the blowing up of $X$ is given as the strict transform

$$\begin{array}{ccc} \tilde{X} & \xrightarrow{h} & X \\ \downarrow & & \downarrow \\ V & \xrightarrow{\pi} & R^n \times R^k \end{array} \tag{1}$$

of $X$ under $\pi$. By [**Sch 4**, 2.2.1] this strict transform $\tilde{X}$ is a closed analytic subset of $V$. Since closed immersions are clearly algebraic and since the composition of two algebraic maps is again such, we only need to prove the claim for the map $\pi$ above.

Without loss of generality we can restrict ourselves to one chart over which this map $\pi$ becomes the map

$$\theta : R^n \times R^k \to R^n \times R^k : (x, t) \mapsto (x_1, x_1 x_2, \ldots, x_1 x_n, t). \tag{2}$$

As can be seen, this is clearly algebraic. ∎

**1.2. Proposition.** *Let $h : \tilde{X} \to X$ be a rigid analytic blowing up map and let $\Sigma$ be a subanalytic subset of $\tilde{X}$ which is closed for the canonical topology. Then $h(\Sigma)$ is also closed for the canonical topology.*

*Proof.* We reduce again to the affinoid case. Using the strict transform diagram (1) in the proof of **(1.1)**, we reduce to the case that $X = R^n \times R^k$ and $0 \times R^k$ is the centre, since closed immersions are clearly closed in the canonical topology. Working on one chart, we further reduce to the case that our map is $\theta$ as given in (2) of above. (Observe that $\operatorname{Im}(\theta)$ is closed for the canonical topology).

So, let $\Sigma$ be a subanalytic subset of $R^n \times R^k$ which is closed for the canonical topology. We need to show that $\theta(\Sigma)$ is closed in $\operatorname{Im}(\theta)$ for the canonical topology. Take any $\alpha \notin \theta(\Sigma)$ but lying in the image of $\theta$. We need to find a small open ball around $\alpha$ not intersecting $\theta(\Sigma)$. There are two possibilities.

*Case 1.* $\alpha \notin 0 \times R^k$. Since $\alpha$ does not lie in the centre, there is a unique element $\beta \in R^n \times R^k$ which maps to $\alpha$ under $\theta$. This element cannot belong to $\Sigma$ neither to $\theta^{-1}(0 \times R^k) = 0 \times R^{n-1} \times R^k$. Since the two latter sets are closed in the canonical topology, we can find a small open ball $B$ around $\beta$ intersecting neither set. Therefore $\theta|_B$ is a homeomorphism onto the open set $\theta(B)$, and the latter contains $\alpha$ and does not intersect $\theta(\Sigma)$. This concludes the first case.

*Case 2.* Assume now that $\alpha$ lies in the centre $0 \times R^k$, say $\alpha = (0, t)$. By assumption, its fiber $\theta^{-1}(\alpha)$ does not intersect $\Sigma$. By **(0.3)**, we therefore can find an $\delta > 0$, such that for any $\tilde{\alpha} \in \theta^{-1}(\alpha)$ and $\tilde{\beta} \in \Sigma$, we have that

$$|\tilde{\alpha} - \tilde{\beta}| \geq \delta. \tag{1}$$

Let $B$ be the open ball $\{\beta \in R^n \times R^k \mid |\beta - \alpha| < \delta\}$ and let $\beta = (b, u) \in B$ be any point. Suppose that $\beta \in \theta(\Sigma)$. Hence there exists an element $\tilde{\beta} = (b_1, b', u) \in \Sigma$, such that $\theta(\tilde{\beta}) = \beta$. Put $\tilde{\alpha} = (0, b', t)$, which lies in the fiber $\theta^{-1}(\alpha)$. It is easy to see that

$$|\tilde{\beta} - \tilde{\alpha}| \leq |\beta - \alpha| < \delta.$$

But this contradicts (1). So $B \cap \theta(\Sigma) = \emptyset$, as we wanted to show. ∎

## 2. Closure of Constructible Sets

**2.1. Lemma.** *Let $A \to B$ be a flat morphism and let $\mathfrak{p}$ be a prime of $A$. Then either $\mathfrak{p}B = B$ or $\mathfrak{p}B \cap A = \mathfrak{p}$.*

*Proof.* Suppose that $\mathfrak{p}B \neq B$. Hence there exists a maximal ideal $\mathfrak{N}$ of $B$, such that $\mathfrak{p}B \subset \mathfrak{N}$. Let $\mathfrak{m} = \mathfrak{N} \cap A$. Then the map $A_\mathfrak{m} \to B_\mathfrak{N}$ is faithfully flat. Hence

$$\mathfrak{p}B_\mathfrak{N} \cap A_\mathfrak{m} = \mathfrak{p}A_\mathfrak{m}.$$

Therefore, if $\alpha \in \mathfrak{p}B \cap A$, then $\alpha \in \mathfrak{p}A_\mathfrak{m}$. In other words, there exists an $s \in A$, with $s \notin \mathfrak{m}$, such that $s\alpha \in \mathfrak{p}$. Since $\mathfrak{p}B \subset \mathfrak{N}$, we must have that $s \notin \mathfrak{p}$ and therefore $\alpha \in \mathfrak{p}$. ∎

**2.2. Corollary.** *Let $X$ be a rigid analytic variety and $V$ and $W$ two closed analytic subsets of $X$, such that no irreducible component of $V$ lies in $W$. Then $\mathrm{cl}(V \setminus W) = V$. In particular, the closure of a constructible set is closed analytic.*

*Remark.* With a *constructible* set, we mean here a finite union of differences of closed analytic subsets.

*Proof.* We easily reduce to the case that $X = \mathrm{Sp}\, A$ is affinoid. Let $x \in V$ and $D = \mathrm{Sp}\, B$ an open ball containing $x$. We need to show that $D \cap (V \setminus W) \neq \emptyset$. Suppose the contrary, so that

$$D \cap V \subset W. \tag{1}$$

Let $\mathfrak{p}$ be the prime ideal of $A$ defining an irreducible component $Z$ of $V$ containing $x$. Let $\mathfrak{a}$ be the radical ideal of $A$ defining $W$. In terms of these ideals, (1) reads as $\mathfrak{a}B \subset \mathrm{rad}(\mathfrak{p}B)$, so that there exists some $n$, such that

$$\mathfrak{a}^n B \subset \mathfrak{p}B. \tag{2}$$

From **(2.1)**, we obtain that $\mathfrak{p}B \cap A = \mathfrak{p}$. (Note that $x \in D \cap Z \neq \emptyset$, implying that $\mathfrak{p}B \neq B$). Hence, from (2) it follows that $\mathfrak{a}^n \subset \mathfrak{a}^n B \cap A \subset \mathfrak{p}$, implying that $Z \subset W$, which contradicts our assumption.

The last assertion now follows immediately from this, since in a difference $V \setminus W$ of closed analytic subsets, we can always discard an irreducible component of $V$ lying in $W$ and hence assume there is none such. ∎

## 3. Closure of Semianalytic Sets

**3.1. Lemma.** *Let $X$ be a rigid analytic variety and $\Sigma \in \mathcal{A}_X^{\mathrm{loc}}$. Then $\Sigma$ is semialgebraic and is clopen for the canonical topology.*

*Proof.* We easily reduce to the case that $X = \operatorname{Sp} A$ is affinoid and that $\Sigma \in \mathcal{A}_X$. Let $\Sigma$ be given as

$$\Sigma = \{\, x \in X \mid |f(x)| \diamond 1 \,\},$$

where $f \in A$ and $\diamond \in \{\leq, <, \geq, >\}$. In order to prove that $\Sigma$ is semialgebraic, we can furthermore assume that $X = R^n$, since any affinoid variety can be viewed as a closed immersion in some $R^n$. Let $f = \sum a_i S^i$, where $S = (S_1, \ldots, S_n)$ and $a_i \in K$. Since the complement of a semialgebraic set (by algebraic Quantifier Elimination) is semialgebraic again, we only have to deal with the two cases that $\diamond$ equals $\leq$ or $<$. Let $\bar{\diamond}$ be defined as $>$ in the first case and as $\geq$ in the second case. Put

$$\tilde{f} = \sum_{|a_i| \bar{\diamond} 1} a_i S^i.$$

By definition of strictly convergent power series, this is a polynomial. We claim that $\Sigma$ coincides with the set

$$\{\, x \in X \mid |\tilde{f}(x)| \diamond 1 \,\}.$$

Put $g = f - \tilde{f}$, then its Gauss norm satisfies $|g| \diamond 1$, so that for each point $x \in X$, we have that $|g(x)| \diamond 1$. From this our claim follows readily, showing that $\Sigma$ is semialgebraic, since $\tilde{f}$ is polynomial.

So it remains to prove that $\Sigma$ is clopen for the canonical topology. Since the complement of a clopen is obviously clopen, we only need to consider the cases that $\diamond$ equals $\leq$ or $\geq$. But then $\Sigma$ is actually an affinoid subdomain whence open in the canonical topology by [**BGR**, 7.2.3. Proposition 1]. We leave it as an exercise to the reader to prove that $\Sigma$ is also closed in the canonical topology. ∎

**3.2. Proposition.** *Let $X$ be a quasi-compact rigid analytic manifold and $\Sigma \in \mathcal{S}_X^{\mathrm{loc}}$. Then there exists a finite admissible affinoid covering $\{U_i\}_i$ of $X$ and maps of rigid analytic varieties $h_i \colon \tilde{X}_i \to U_i$, which are compositions of finitely many blowing up maps, such that, for each $i$, we can find a set $A_i \in \mathcal{A}_{\tilde{X}_i}^{\mathrm{loc}}$ and a closed analytic subset $V_i \subset \tilde{X}_i$, so that we can write $h_i^{-1}(\Sigma)$ as either $A_i \cup V_i$ or $A_i \setminus V_i$.*

*Proof.* By quasi-compactness, we easily reduce to the case that $X = \operatorname{Sp} A$ is affinoid and that $\Sigma \in \mathcal{S}_X$. Let $\Sigma$ be given as
$$\Sigma = \{\, x \in X \mid |f(x)| \diamond |g(x)| \,\},$$
where $f, g \in A$ and $\diamond \in \{\leq, <, \geq, >\}$. By [**Sch 2**, Lemma 4.2] applied to the couple $(f, g)$, we can find an admissible affinoid covering $\{U_i\}_i$ of $X$, compositions of blowing up maps $h_i \colon \tilde{X}_i \to U_i$ and finite admissible affinoid coverings $\{\tilde{Y}_{ij}\}_j$ of $\tilde{X}_i$, such that, for each $i$ and each $j$, either $(g \circ h_i)|_{\tilde{Y}_{ij}}$ divides $(f \circ h_i)|_{\tilde{Y}_{ij}}$, or vice versa, $(f \circ h_i)|_{\tilde{Y}_{ij}}$ divides $(g \circ h_i)|_{\tilde{Y}_{ij}}$ as elements in $\mathcal{O}_{\tilde{X}_i}(\tilde{Y}_{ij})$.

Fix an $i$ and let us write $h \colon \tilde{X} \to U$ for $h_i \colon \tilde{X}_i \to U_i$ and $\tilde{Y}_j$ for $\tilde{Y}_{ij}$. Renumber the $\tilde{Y}_j$ in such way that $g \circ h$ divides $f \circ h$ on $\tilde{Y}_j$ for $j = 1, \ldots, s$, and $f \circ h$ divides $g \circ h$ on $\tilde{Y}_j$ for $j = s+1, \ldots, t$, where we do not indicate restrictions anymore in order not to overload our notation.

*Case 1.* $1 \leq j \leq s$. So there exists a $\nu_j \in \mathcal{O}_{\tilde{X}}(\tilde{Y}_j)$ such that $f \circ h = (g \circ h) \cdot \nu_j$ on $\tilde{Y}_j$ Hence we have that

(1) $$h^{-1}(\Sigma) \cap \tilde{Y}_j = \left\{\, x \in \tilde{Y}_j \mid |g(h(x)) \cdot \nu_j(x)| \diamond |g(h(x))| \,\right\}.$$

We will refer to the cases $\diamond$ equals $\geq$ or $\leq$ as the *open* cases. Therefore, in the open cases we calculate from (1) that

(Open) $$h^{-1}(\Sigma) \cap \tilde{Y}_j = \left\{\, x \in \tilde{Y}_j \mid |\nu_j(x)| \diamond 1 \,\right\} \cup \mathrm{V}(f \circ h, g \circ h).$$

The remaining cases that $\diamond$ equals $<$ or $>$, will be called the *closed* cases and we calculate that then

(Closed) $$h^{-1}(\Sigma) \cap \tilde{Y}_j = \left\{ x \in \tilde{Y}_j \mid |\nu_j(x)| \diamond 1 \right\} \setminus \mathrm{V}(f \circ h, g \circ h).$$

*Case 2.* $s+1 \leq j \leq t$. Hence there exists a $\nu_j \in \mathcal{O}_{\tilde{X}}(\tilde{Y}_j)$ such that $g \circ h = (f \circ h) \cdot \nu_j$ on $\tilde{Y}_j$. We obtain, by symmetry, exactly the same description (Open) respectively (Closed), for the open, respectively the closed cases, with the only difference that now $\nu_j$ and $1$ are interchanged. It is not too hard to check, that on intersections $\tilde{Y}_j \cap \tilde{Y}_k$, the functions $\nu_j$ and $\nu_k$ agree on these intersections, if $j$ and $k$ are both at most $s$ or both strictly bigger than $s$. They are each others inverse in the remaining case. From this it follows that we can construct a set $A \in \mathcal{A}_{\tilde{X}}^{\mathrm{loc}}$, such that $h^{-1}(\Sigma)$ equals $A \cup \mathrm{V}(f \circ h, g \circ h)$ in the open case, or $A \setminus \mathrm{V}(f \circ h, g \circ h)$ in the closed case. ∎

*Remark.* If we would use a generalisation of [**Sch 2**, Lemma 4.2] in which finitely many couples $f_i, g_i$, rather than just one such couple, are brought simultaneously into relative division after blowing up, then we can improve our proposition to find the claimed admissible covering and maps so that they work simultaneously for finitely many sets $\Sigma_j \in \mathcal{S}_X^{\mathrm{loc}}$. To prove this generalised version of [**Sch 2**, Lemma 4.2], just take for $F$ in the proof of loc. cit., the product of all $f_i g_i (f_i - g_i)$ instead of just $fg(f-g)$ and the proof carries over almost verbatim.

**3.3. Theorem.** *Let $X$ be a rigid analytic variety and $\Sigma$ a semianalytic subset of $X$. Then its closure $\mathrm{cl}(\Sigma)$ in the canonical topology is also semianalytic.*

*Proof.* We can again reduce to the case that $X$ is affinoid and that $\Sigma$ is globally semianalytic in $X$. After taking a closed immersion, we can even assume that $X = R^n$, for some $n$. Since the closure of a finite union of sets is the union of the closure of these sets, we can even assume that $\Sigma$ is basic. In other words, we can write $\Sigma$ as a finite intersection

$$\Sigma = \bigcap_{j=1}^{s} \Sigma_j,$$

where $\Sigma_j \in \mathcal{S}_{R^n}$. By the remark following **(3.2)**, we can find a finite admissible affinoid covering $\{U_i\}_i$ of $R^n$ and maps of rigid analytic varieties $h_i \colon \tilde{X}_i \to U_i$ which are compositions of finitely many blowing up maps, such that each $h_i^{-1}(\Sigma_j)$ is either of the form $A \cup V$ or $A \setminus V$, with $A \in \mathcal{A}_{\tilde{X}_i}^{\mathrm{loc}}$ and $V \subset \tilde{X}_i$ a closed analytic subset. Dropping indices, let $h \colon \tilde{X} \to U$ be one of these maps $h_i$. Renumber the $\Sigma_j$ in such way that,

$$\begin{aligned} h^{-1}(\Sigma_j) &= A_j \cup V_j & j &= 1, \ldots, t \\ &= A_j \setminus V_j & j &= t+1, \ldots, s. \end{aligned}$$

where the $A_j \in \mathcal{A}_{\tilde{X}}^{\mathrm{loc}}$ and where the $V_j \subset \tilde{X}$ are closed analytic subsets. In other words, we have that

(1) $$h^{-1}(\Sigma) = (A_1 \cup V_1) \cap \cdots \cap (A_t \cup V_t) \cap (A_{t+1} \setminus V_{t+1}) \cap \cdots \cap (A_s \setminus V_s).$$

Applying distrubitivity of $\cap$ over $\cup$, we can rewrite this as

(2) $$h^{-1}(\Sigma) = \bigcup_{k=1}^{q} B_k \cap (W_k \setminus Z_k),$$

where the $B_k$ are finite intersections of sets in $\mathcal{A}_{\tilde{X}}^{\mathrm{loc}}$ and where $W_k$ and $Z_k$ are closed analytic subsets of $\tilde{X}$. Without loss of generality, we can even assume, for each $k = 1, \ldots, q$, that no irreducible component of $W_k$ lies in $Z_k$. Using (0.1) and **(3.1)**, we can now calculate from (2) the closure of $h^{-1}(\Sigma)$ in the canonical topology as

(3) $$\mathrm{cl}(h^{-1}(\Sigma)) = \bigcup_{k=1}^{q} B_k \cap \mathrm{cl}(W_k \setminus Z_k).$$

From **(2.2)** we learn that $\operatorname{cl}(W_k \setminus Z_k) = W_k$. From **(1.2)** it follows that each $h(B_k \cap W_k)$ is closed for the canonical topology in $U$. Moreover, since the image of a semialgebraic set under a semialgebraic map is semialgebraic in the image of the map and since both the closed immersion $W_k \hookrightarrow \tilde{X}$ and $h$ are algebraic by **(2.1)** and whence their composition, we conclude that $h(B_k \cap W_k)$ is semialgebraic in $h(W_k)$ by **(3.1)**. Since $h$ is a composition of proper maps by [**Sch 4**, Theorem 3.2.1], $h(W_k)$ is a closed analytic subset of $U$ by [**BGR**, 9.6.3. Proposition 3]. From this it follows that $h(B_k \cap W_k)$ is semianalytic in $U$. Just observe that semianalytic (and hence semialgebraic) sets are sent to semianalytic sets under closed immersions, see [**Sch 2**, Proposition 1.3]. Hence both closedness as well as semianalyticity also hold for their union, which equals $h(\operatorname{cl}(h^{-1}(\Sigma)))$ by (3).

Therefore, using (0.2), we have inclusions

$$(4) \qquad \Sigma \cap U \subset h(\operatorname{cl}(h^{-1}(\Sigma))) \subset \operatorname{cl}(\Sigma) \cap U.$$

Since the middle set of (4) is also closed for the canonical topology, we conclude that the two last sets of (4) must coincide, proving that $\operatorname{cl}(\Sigma) \cap U$ is semianalytic. Returning back to our original covering, i.e., writing back subscripts, we found that each $\operatorname{cl}(\Sigma) \cap U_i$ is semianalytic in $U_i$, proving that $\operatorname{cl}(\Sigma)$ is semianalytic in $R^n$. ∎

## 4. BOUNDARIES

**4.1. Lemma.** *Let $X = \operatorname{Sp} A$ be an affinoid variety and let $f, g \in A$. Suppose there exists an admissible affinoid covering $\{X_i\}_i$ of $X$ and surjective maps of rigid analytic varieties $h_i \colon Y_i \twoheadrightarrow X_i$, such that for each $i$ we have that $f \circ h_i$ divides $g \circ h_i$ on $Y_i$. Then $\operatorname{V}(f) \subset \operatorname{V}(g)$.*

*Proof.* Let $x \in \operatorname{V}(f)$. By assumption, there exists an $i$, such that $x \in h_i(X_i)$. Choose $y \in Y_i$, such that $h_i(y) = x$. Since $f \circ h_i$ divides $g \circ h_i$ and since $f(x) = f(h_i(y)) = 0$, we also have that $g(x) = 0$, as we needed to show. ∎

**4.2. Proposition.** *Let $X = \operatorname{Sp} A$ be an affinoid manifold and let $f, g \in A$. Let*

$$\Sigma = \{\, y \in X \mid |f(y)| < |g(y)| \,\}.$$

*Let $x$ be a point in $\operatorname{V}(f, g)$ which is not contained in the closure $\operatorname{cl}(\Sigma)$ of $\Sigma$ in the canonical topology. Then there exists an admissible affinoid $U$ of $X$, containing $x$, such that*

$$\operatorname{V}(f, g) \cap U = \operatorname{V}(f) \cap U.$$

*Proof.* Let $F = fg(f - g)$. By Embedded Resolution of Singularities [**Sch 5**, Theorem 3.2.3], we can find a finite admissible affinoid covering $\{W_i\}_{i=1,\ldots,s}$ of $X$ and (surjective) maps $h_i \colon \tilde{X}_i \to W_i$ of rigid analytic manifolds, where each $h_i$ is a finite composition of blowing up maps, such that $F \circ h_i$ has normal crossings in $\tilde{X}_i$. From the proof of [**Sch 2**, Lemma 4.2] it then follows that there exists a finite admissible affinoid covering $\{\tilde{Y}_{ij}\}_j$ of $\tilde{X}_i$, such that either $f \circ h_i$ divides $g \circ h_i$ on $\tilde{Y}_{ij}$ or vice versa.

We claim that we can find, for each $i = 1, \ldots, s$, an admissible open $U_i$ of $X$, containing $x$, such that $f \circ h_i$ divides $g \circ h_i$ on the whole of $h_i^{-1}(U_i)$. Suppose this has been proven. Let $U$ be the intersection of all $U_i$, which is again an admissible affinoid containing $x$. We can now apply **(4.1)** to $U$ and the collection of restrictions $h_i|_{h_i^{-1}(U)}$ to conclude the proof.

So we only need to show the claim. Fix an $i = 1, \ldots, s$ and let us write $h \colon \tilde{X} \to W$ for $h_i \colon \tilde{X}_i \to W_i$ and $\tilde{Y}_j$ for $\tilde{Y}_{ij}$. The only non-trivial case is when $x \in W$ and $g \circ h$ divides $f \circ h$ on $\tilde{Y}_j$. Let $\tilde{x} \in \tilde{Y}_j$ be such that $h(\tilde{x}) = x$. There exists a $\nu \in \mathcal{O}_{\tilde{X}}(\tilde{Y}_j)$, such that $f \circ h = (g \circ h) \cdot \nu$. By definition of normal crossings, there exists a regular system of parameters $\xi = (\xi_1, \ldots, \xi_s)$ in $\mathcal{O}_{\tilde{X}, \tilde{x}}$, such that

$$
\begin{aligned}
f \circ h &= u \xi^\alpha \\
g \circ h &= u \xi^\beta \\
\nu &= \xi^\Omega,
\end{aligned}
\tag{1}
$$

where $u$ is a unit in $\mathcal{O}_{\tilde{X},\tilde{x}}$ and $\alpha$, $\beta$ and $\Omega$ are indices with $\alpha = \beta + \Omega$. Hence there exists a Zariski open $D_{\tilde{x}}$ in $\tilde{Y}_{\tilde{j}}$ containing $\tilde{x}$, such that $\nu$ and $u$ are defined over $D_{\tilde{x}}$ and such that the equations (1) still hold in $\mathcal{O}_{\tilde{X}}(D_{\tilde{x}})$, with $u$ a unit in the latter ring.

Since $x \notin \text{cl}(\Sigma)$, it follows from (0.2) that $\tilde{x}$ does not belong to $\text{cl}(h^{-1}(\Sigma))$. Let us calculate $\text{cl}(h^{-1}(\Sigma))$. It is easy to see that

$$h^{-1}(\Sigma) \cap \tilde{Y}_{\tilde{j}} = \left\{ y \in \tilde{Y}_{\tilde{j}} \mid |\nu(y)| < 1 \wedge g(h(y)) \neq 0 \right\}.$$

From this it is not too hard to deduce that

(2) $$\text{cl}(h^{-1}(\Sigma)) \cap \tilde{Y}_{\tilde{j}} = \left\{ y \in \tilde{Y}_{\tilde{j}} \mid |\nu(y)| < 1 \right\}.$$

From the fact that $\tilde{x} \notin \text{cl}(h^{-1}(\Sigma))$, we obtain by (2) that $1 \leq |\nu(\tilde{x})|$. On the other hand, we have by (1), that $\nu = \xi^{\Omega}$. But each $\xi_i$ vanishes in the point $\tilde{x}$. Therefore, we must have that $\Omega = 0$, so that actually $f \circ h$ divides $g \circ h$, when restricted to $D_{\tilde{x}}$.

Let $D$ be the union of all the $D_{\tilde{x}}$ when $\tilde{x}$ runs over the fiber $h^{-1}(x)$. Hence $D$ is again Zariski open in $\tilde{X}$ and $f \circ h$ divides $g \circ h$ on $D$. In particular, the complement $V$ of $D$ is a closed analytic subset. Since $h$ is the composition of proper maps by [**Sch 4**, Theorem 3.2.1], it follows from [**BGR**, 9.6.3. Proposition 3] that $h(D)$ is closed analytic in $W$ and hence in particular closed for the canonical topology. Since $h^{-1}(x) \subset D$, it follows that $x \notin h(V)$. Therefore, there exists an admissible affinoid $U$, containing $x$, such that $U \cap h(V) = \emptyset$. But this clearly implies that $h^{-1}(U) \subset D$, so that $U$ is the required admissible affinoid.  ∎

**4.3. Lemma.** *Let $X = \text{Sp}\,A$ be an affinoid variety and $U = \text{Sp}\,C$ an admissible affinoid in $X$. Suppose that both $A$ and $C$ are unique factorisation domains. Let $f, g \in A$ be relative prime in $A$, then the same holds for their restrictions $f|_U$ and $g|_U$ in $C$.*

*Remark.* Recall that in an unique factorisation domain $B$, any two elements $r$ and $s$ have a greatest common divisor (see for instance [**Mats**, Theorem 20.5]). If one and hence all greatest common divisors of $r$ and $s$ are units, then we say that $r$ and $s$ are *relatively prime* or *have no common factor*.

*Proof.* Suppose that $f$ and $g$ have a greatest common divisor $d \in C$ which is not a unit, (where we have written $f$ for $f|_U$ and likewise for $g$). Hence there exist $\tilde{f}, \tilde{g} \in C$ such that $f = \tilde{f} \cdot d$ and $g = \tilde{g} \cdot d$ in $C$. In particular, we have that

$$\tilde{g} \cdot f - \tilde{f} \cdot g = 0$$

in $C$. Since $C$ is flat over $A$ (see [**BGR**, 7.3.2. Corollary 6]), we can find, for $i = 1, \ldots, s$, elements $f_i, g_i \in A$ and $r_i, s_i \in C$, such that

(1) $$g_i \cdot f - f_i \cdot g = 0$$

in $A$ and

(2) $$\tilde{f} = \sum_{i=1}^{s} r_i f_i$$
$$\tilde{g} = \sum_{i=1}^{s} s_i g_i$$

in $C$.

Since $f$ and $g$ are relatively prime in $A$, it follows from (1), that each $f_i$ is divisible by $f$ in $A$. Hence, by (2), $\tilde{f}$ would be divisible by $f$ in $C$, implying that $d$ has to be a unit, contradiction.  ∎

**4.4. Theorem.** *Let $X$ be a rigid analytic variety and let $\Sigma \in \mathcal{S}_X^{\mathrm{loc}}$. Then the boundary of $\Sigma$ is the union of a closed analytic subset and the intersection of a closed analytic hypersurface with a (semianalytic) open set.*

*More precisely, suppose*
$$\Sigma = \{\, x \in R^n \mid |f(x)| < |g(x)| \,\},$$
*where $f, g \in K\langle S\rangle$ are relatively prime. Then the boundary of $\Sigma$ is the analytic subset $\mathrm{V}(f,g)$.*

*Proof.* As before we can assume already that $X = R^n$, for some $n$, and that $\Sigma \in \mathcal{S}_{R^n}$. By considering complements, we only need to treat the case that
$$\Sigma = \{\, x \in R^n \mid |f(x)| < |g(x)| \,\},$$
where $f, g \in K\langle S\rangle$. In particular, one checks immediately that $\Sigma$ is open in the canonical topology and $\Sigma \cup \mathrm{V}(f,g)$ is closed in the canonical topology. In case $f$ is a unit in $K\langle S\rangle$, we have that $\Sigma \in \mathcal{A}_X$, whence a clopen in the canonical topology by **(3.1)**. So we are done in this case.

Therefore, assume that $f$ is not a unit. Let us first consider the case that $f$ and $g$ are relatively prime. We have proven the theorem for this case if we show that

(1) $$\mathrm{cl}(\Sigma) = \Sigma \cup \mathrm{V}(f,g),$$

since $\mathrm{V}(f,g)$ would then be the boundary of $\Sigma$. We already observed that the right hand side of (1) is closed for the canonical topology, which establishes the $\subset$-inclusion. Therefore, assume that (1) does not hold, so that there exists an element $y \in \mathrm{V}(f,g) \setminus \mathrm{cl}(\Sigma)$. Hence, by Proposition **(4.2)**, there exists an admissible affinoid $U = \mathrm{Sp}\, C$ of $R^n$, containing $y$, such that
$$\mathrm{V}(f,g) \cap U = \mathrm{V}(f) \cap U.$$

By taking a polydisk contained in $U$, we can even assume that $C$ is a unique factorisation domain. By the Nullstellensatz for $U$ (see for instance [**BGR**, 7.1.2. Theorem 3]), this means that there exists an $n$, such that $g^n \in fC$. But this would imply that $f$ and $g$ are not relatively prime in $C$, since $f$ is not a unit, contradicting our assumption by means of **(4.3)**.

For the general case, we can write $f = hf_0$ and $g = hg_0$ with $f_0, g_0, h \in K\langle S\rangle$ such that $f_0$ and $g_0$ are relatively prime and $h \neq 0$. Put
$$\Sigma_0 = \{\, x \in R^n \mid |f_0(x)| < |g_0(x)| \,\}.$$

The reader can check that then $\Sigma = \Sigma_0 \setminus \mathrm{V}(h)$. In fact, we have a disjoint union

(2) $$\Sigma_0 = \Sigma \cup (\Sigma_0 \cap \mathrm{V}(h)).$$

We claim that

(3) $$\mathrm{cl}(\Sigma) = \mathrm{cl}(\Sigma_0).$$

Indeed, suppose $x \in \mathrm{cl}(\Sigma_0)$ so that there exists an open $U$ of $X$ containing $x$ with $U \cap \Sigma_0$ non-empty. This latter set is also open, since $\Sigma_0$ is, and hence $h$ cannot be identical zero on it. In other words, there exists a point $u \in U \cap \Sigma_0$ with $h(u) \neq 0$ and hence $u \in \Sigma$. This finishes the proof of our claim (3). Now, by the relative prime case above we know that $\mathrm{cl}(\Sigma_0) = \Sigma_0 \cup V$, where $V = \mathrm{V}(f_0, g_0)$. Putting this together with (2) and (3) yields
$$\partial \Sigma = V \cup (\Sigma_0 \cap \mathrm{V}(h)),$$
as wanted. ∎

MATHEMATICAL INSTITUTE
UNIVERSITY OF OXFORD
24-29 ST. GILES
OXFORD OX1 3LB (UNITED KINGDOM)
  *Current address*:
Fields Institute
222 College Street
Toronto, Ontario, M5T 3J1 (Canada)
  *E-mail address*: hschoute@fields.utoronto.ca